\documentclass[12pt,a4paper]{article}
\usepackage{amsfonts}
\usepackage{amssymb}
\usepackage{latexsym}

\newcommand{\mc}[1]{\mathcal{#1}}
\newcommand{\mbb}[1]{\mathbb{#1}}

\newcommand{\Hom}{\mathrm{Hom}}
\newcommand{\id}{\mathrm{id}}

\newcommand{\im}{\mathrm{im}}

\newcommand{\ovr}[1]{\mathop{\vbox{\ialign{
				##\crcr
    ${\scriptstyle\hfil\;\;#1\;\;\hfil}$\crcr
    \noalign{\kern-1pt\nointerlineskip}
    \rightarrowfill\crcr}}\;}}
\newcommand{\map}{\longrightarrow}
\newcommand{\mb}[1]{\hbox{$#1$}}

\newcommand{\qmod}[2]{{\hbox{}^{\displaystyle{#1}}}\!\big/\!
\hbox{}_{\displaystyle{#2}}}

\newcommand{\sett}[2]{\{\ #1\ \vert\ #2\ \}}

\newtheorem{thm}{Theorem}%[subsection]

\font\caps=cmb10

\def\TSkip{\bigskip}

\def\Abstract{\begingroup\narrower
    \parskip=\medskipamount\parindent=0pt{\caps Abstract. }}
\def\EndAbstract{\par\endgroup\TSkip\TSkip}

\begin{document}

\title{The natural metric in the Horrocks-Mumford bundle is not Hermitian-Einstein.}         
\author{O.F.B. van Koert, M. L\"ubke}        
\date{}          
\maketitle

\Abstract
{\footnotesize The Horrocks-Mumford bundle $E$ is a famous stable complex vector bundle of rank 2 on 4-dimensional complex projective space. By construction, $E$ has a natural Hermitian metric $h_1$. On the other hand, stability implies the existence of a Hermitian-Einstein metric in $E$ which is unique up to a positive scalar. Now the obvious question is if $h_1$ is in fact the Hermitian-Einstein metric. In this note we indicate how to show by computation that this is not the case.} \EndAbstract

\section{Introduction and main result}

Let $N$ be a null correlation bundle on complex projective 3-space $\mbb{P}_3$, i.e. a quotient of $\Omega^1_{\mbb{P}_3}$ by $\mc{O}_{\mbb{P}_3}(-1)$ (see e.g. \cite{OSS}). Then $N$ is stable in the sense of \cite{OSS}, or equivalently, $g$-stable in the sense of \cite{LT}, where $g$ is the Fubini-Study metric in $\mbb{P}_3$, so the Kobayashi-Hitchin correspondence tells us that there exists a $g$-Hermitian-Einstein metric $h_0$ in $N$, which is unique up to a constant positive factor. On the other hand, the standard metric in $\mbb{C}^4$ not only induces the Fubini-Study metric in $\mbb{P}_3$, but also natural metrics in $\Omega^1_{\mbb{P}_3}$ and $\mc{O}_{\mbb{P}_3}(-1)$, and hence a metric $h_1$ in the quotient $N$, too. Now the obvious question arises: 

{\bf (Q)} Does it hold\ \ \mb{h_1 = c\cdot h_0}\ \ with a positive constant $c$, or equivalently, does $h_1$ satisfy the $g$-Hermitian-Einstein equation
\[
(HE)\ \ \ K_{h_1} = \lambda\cdot \id_E\ 
\]
where $K_{h_1}$ is the mean curvature of $h_1$ and $\lambda$ a real constant? 

This question was answered in the affirmative in \cite{L} by manual computations with respect to local coordinates and a local holomorphic frame field.

In this note we consider the following similar situation. It is well known that on the 4-dimensional complex projective space\ \ \mb{\mbb{P}_4 = \mbb{P}(\mbb{C}^5)}\ \ there exists a stable holomorphic 2-bundle $E$ with Chern numbers\ \ \mb{c_1(E) = 5}\ \ and\ \ \mb{c_2(E) = 10\ ,}\ the {\it Horrocks-Mumford bundle} \cite{HM}, \cite{OSS}.  Again, stability of $E$ in the sense of \cite{OSS} is the same as $g$-stability in the sense of \cite{LT}, where $g$ is the Fubini-Study metric in $\mbb{P}_4$, so there exists a $g$-Hermitian-Einstein metric $h_0$ in $E$, which is unique up to a constant positive factor. 
On the other hand, using the construction of $E$ given in \cite{OSS} one gets in a natural way an explicit metric $h_1$ in $E$ induced by the standard metric in $\mbb{C}^5$, hence question (Q) arises for $E$, too. Again we used explicit calculations in local coordinates to tackle this problem, and the result is 

\begin{thm}
The metric natural metric $h_1$ in the Horrocks-Mumford bundle is {\bf NOT} $g$-Hermitian-Einstein.
\end{thm}

In section 2 we sketch our approach to the problem, and in section 3 we give some details and explicit formulae which should be sufficient to make our calculations reproducible.

\section{Our approach} 

The construction in [OSS] we use does not produce the bundle $E$ directly, but the bundle\ \ \mb{E(-2) = E\otimes\mc{O}_{\mbb{P}_4}(-2)} and a metric $h$ in it. Let $h_2$ be the standard metric in\ \ \mb{\mc{O}_{\mbb{P}_4}(2) = \mc{O}_{\mbb{P}_4}(1)^{\otimes 2}\ ,}\ induced by the canonical inclusion\ \ \mb{\mc{O}_{\mbb{P}_4}(1)^* = \mc{O}_{\mbb{P}_4}(-1) \hookrightarrow \mbb{P}_4\times\mbb{C}^5}\ \ and the standard metric in $\mbb{C}^5$, then the natural metric $h_1$ in\ \ \mb{E = E(-2)\otimes\mc{O}_{\mbb{P}_4}(2)}\ \ is the metric induced by $h$ and $h_2$. 

Our initial guess was that $h_1$ would be indeed $g$-Hermitian-Einstein. Since $h_2$ is known to be $g$-Hermitian-Einstein, this is equivalent to $h$ being $g$-Hermitian-Einstein. Hence we attempted to show that the equation $(HE)$ holds for $h$; by continuity, it suffices to do that on some open dense subset $U_0^*$ of $\mbb{P}_4$. So in a suitable chart for $\mbb{P}_4$ we explicitely determined a matrix representation $H$ of the metric $h$ with respect to a holomorphic frame field; this already involved algebraic calculations which where impossible to do by hand ($H$ contains rational expressions in the 8 real variables $x_1,\bar{x}_1,\ldots,x_4,\bar{x}_4$, with numerators of degree up to 16), so we used the computer package MAPLE. The next step would have been the calculation of the mean curvature, i.e. essentially the matrix\ \ \mb{K = (K_{ij})_{i,j=1}^2}\ \ where
\[
K_{ij} = -\sum\limits_{\alpha,\beta = 1}^{4}g^{\beta\alpha}\left(\sum\limits_{k = 1}^2\frac{\partial^2H_{ik}}{\partial x_\alpha\partial \bar{x}_\beta}H^{kj} - \sum\limits_{k,l,m = 1}^2\frac{\partial H_{ik}}{\partial x_\alpha}H^{kl}\frac{\partial H_{lm}}{\partial \bar{x}_\beta}H^{mj}\right)\ .
\]
Here the $g_{\alpha\beta}$ are the coefficients of the Fubini-Study metric with respect to the local holomorphic coordinates $x_\alpha$, and upper indices mean coefficients of the inverse matrix. Now if the metric $h$ was $g$-Hermitian, then since\ \ \mb{c_1(E(-2)) = 1}\ \ this would be equivalent to\ \ \mb{K(x) = \left(\begin{array}{cc}2&0\\0&2\end{array}\right)}\ \ for all\ \ \mb{x \in U_0^*\ .}\ (Notice that the constant $\lambda$ in $(HE)$ is determined by the topology of $E$ (see e.g. [LT]) and can therefore be determined a priori.)
 
Unfortunately, MAPLE was not able (at least on our computer) to calculate $K$ in a general point $x$ (the main problem being the inverse of $H$), so we decided to do some testing. 

For this, we first let MAPLE determine the derivatives involved in the formula for $K_{ij}$ in a general point. Then we took the particular point\ \ \mb{x_0 = (x_1^0,\ldots,x_4^0) = (1,1,1,1)\ ,}\ and could calculate $K(x^0)$ (inversion of the scalar matrix $H(x^0)$ is easy). The result was (as we had hoped) indeed
\[
K(x^0) = \left(\begin{array}{cc}2&0\\0&2\end{array}\right)\ .
\]
We repeated the procedure with the point $x^1 = (2,1,1,1)$, and again we got
\[
K(x^1) = \left(\begin{array}{cc}2&0\\0&2\end{array}\right)\ .
\]
This seemed to indicate that we had in fact some chance to be right with our first guess.

In the meantime, we had started a different test by looking at the determinant line bundle\ \ \mb{L := \det E(-2)\ .}\ The induced metric $\det h$ in $L$ is given over $U_0^*$ by the function $\det H$, and if $h$ was $g$-Hermitian-Einstein, then $\det h$ would be $g$-Hermitian-Einstein, too; more precisely, the mean curvature $K_{\det H}$ of $\det h$ would be the constant function 4. Again we got a problem: MAPLE could calculate $\det H$, but was not able to simplify the resulting rational function to a form from which it could determine $K_{\det H}$. But we made the motivated guess that
\[
\det H = \frac{|x_1|^2|x_2|^2|x_3|^2|x_4|^2}{1+\|x\|^2}\ ,
\]
and where able (using MAPLE) to verify the correctness of this formula. Now it was easy to check (even by hand) that indeed\ \ \mb{K_{\det H} \equiv 4\ ,}\ i.e. that the induced metric in $\det E(-2)$, and hence that in $\det E$, is $g$-Hermitian-Einstein.
 
So far everything seemed to be okay, but testing of a third point\ \ \mb{x^2 = (1+i,1,1,1)}\ \ gave the disappointing result
\[
K(x^2) = \left(\begin{array}{cc}2&-\frac{217}{25992}\\-\frac{217}{25992}&2\end{array}\right)\ !
\]
This of course meant precisely what we did \underbar{not} want to show, namely that the metric in $E(-2)$, and hence the metric in $E$, is not $g$-Hermitian-Einstein.

\section{Some details and formulae}

The bundle $E(-2)$ is the cohomology of a monad
\[
(M)\ \ \ 0 \rightarrow \mbb{C}^5\otimes\mc{O}_{\mbb{P}_4} \ovr{a} (\Lambda^2Q)^{\oplus 2} \ovr{b} (\mbb{C}^5)^*\otimes\Lambda^4Q \rightarrow 0\ ,
\]
i.e. $a$ is injective, $b$ is surjective, it holds\ \ \mb{\im(a) \subset \ker(b)\ ,}\ and\ \mb{E(-2) = \qmod{\ker(b)}{\im(a)}\ ;}\ here\ \ \mb{Q = T_{\mbb{P}_4}(-1)}\ \ (see [OSS]). Let\ \ \mb{\pi : \mbb{C}^5\otimes\mc{O}_{\mbb{P}_4} \rightarrow Q}\ \ be the natural projection in the Euler sequence. The standard Hermitian inner product in $\mbb{C}^5$ defines the standard flat Hermitian metric in the trivial bundle $\mbb{C}^5\otimes\mc{O}_{\mbb{P}_4}$, and hence a quotient metric $h_Q$ in $Q$. This induces a metric $\Lambda^2h_Q$ in $\Lambda^2Q$, and hence a metric $h_3$ in $(\Lambda^2Q)^{\oplus 2}$ by taking the two summands as orthogonal. Next we get a metric $h_b$ in $\ker(b)$ by restricting $h_3$, and finally a quotient metric $h$ in $E(-2)$.\\
Let\ \ \mb{x = (x_0:x_1:\ldots:x_4)}\ \ be the homogeneous coordinates in $\mbb{P}_4$ with respect to the standard basis $e_0$,\ldots,$e_4$ of $\mbb{C}^5$. The holomorphic section in $\mbb{C}^5\otimes\mc{O}_{\mbb{P}_4}$ defined by $e_i$ is denoted $\tilde{e}_i$, and we define\ \ \mb{v_i := \pi(\tilde{e}_i) \in H^0(\mbb{P}_4,Q)\ ,}\ \mb{i = 0,\ldots,4\ .}\\ 
Over\ \ \mb{U_0 := \sett{x \in \mbb{P}_4}{x_0 \not= 0}\ ,}\ \mb{\underline{v} := (v_1,\ldots,v_4)}\ \ is a holomorphic frame field for $Q$. For\ \ \mb{x = (1:x_1:\ldots:x_4) \in U_0} the quotient metric in $Q(x)$ is given by
\[
h_Q(v_i,v_j)(x) = \delta_{ij}- \frac{\bar{x}_ix_j}{n}\ \ ,\ \ 1 \leq i,j \leq 4\ ,
\] 
where\ \ \mb{n = 1 + \sum\limits_{i=1}^4\vert x_i\vert^2\ .}\ A holomorphic frame field\ \ \mb{\underline{u} = (u_1,\ldots,u_6)}\ \ for $\Lambda^2Q$ over $U_0$ is given by\ \ \mb{u_1 := v_1\wedge v_2\ ,}\ \ \mb{u_2 := v_1\wedge v_3\ ,}\ \ \mb{u_3 := v_1\wedge v_4\ ,}\ \ \mb{u_4 := v_2\wedge v_3\ ,}\ \ \mb{u_5 := v_2\wedge v_4\ ,}\ \ \mb{u_6 := v_3\wedge v_4\ .}\ Since\ \ \mb{\Lambda^2h_Q(v_i\wedge v_j,v_k\wedge v_l) = h_Q(v_i,v_k)h_Q(v_j,v_l) - h_Q(v_i,v_l)h_Q(v_j,v_k)\ ,}\ it is easy to determine the matrix representation of $\Lambda^2h_Q(x)$ with respect to $u(x)$.
The holomorphic frame field\ \ \mb{\underline{b}  := (b_1,\ldots,b_{12})}\ \ for $(\Lambda^2Q)^{\oplus 2}$ is defined by\ \ \mb{b_i := (u_i,0)}\ \ for\ \ \mb{1 \leq i \leq 6}\ \ and\ \ \mb{b_i := (0,u_{i-6})}\ \ for\ \ \mb{7 \leq i \leq 12\ ;}\ then the matrix representation of $h_3(x)$ with respect to $\underline{b}(x)$ is\ \ \mb{h_3(x) = \left(\begin{array}{cc}\Lambda^2h_Q(x) & 0 \\ 0 & \Lambda^2h_Q(x)\end{array}\right)\ .}

Define\ \ \mb{a_\pm : \mbb{C}^5 \map \Lambda^2\mbb{C}^5}\ \ by\ \ \mb{a_+(e_i) := e_{i+2}\wedge e_{i+3}\ ,}\ \mb{a_-(e_i) := e_{i+1}\wedge e_{i+4}\ ,}\ \mb{0 \leq i \leq 4}\ \ (indices $mod$ 5). Then the map $a$ in $(M)$ is defined as the composition
\[
a(x) : \mbb{C}^5 \ovr{(a_+,a_-)} (\Lambda^2\mbb{C}^5)^{\oplus 2} \ovr{(\Lambda^2\pi(x))^{\oplus 2}} (\Lambda^2Q(x))^{\oplus 2}\ .
\]
Using\ \ \mb{\pi(x)(e_0) = -\sum\limits_{i=1}^4 x_iv_i\ ,}\ it follows that the basis $(a_3,\ldots,a_7)$,\ \mb{a_i := a(e_i)\ ,}\ \ of $\im(a)$ is given in coordinates with respect to $\underline b$ as
\begin{eqnarray*}
a_3(x) &=& (0,0,0,1,0,0,0,0,1,0,0,0)\ , \\
a_4(x) &=& (0,0,0,0,0,1,x_1,0,0,-x_3,-x_4,0)\ , \\
a_5(x) &=& (0,0,x_1,0,x_2,x_3,0,-1,0,0,0,0)\ , \\
a_6(x) &=& (x_2,x_3,x_4,0,0,0,0,0,0,0,-1,0)\ , \\
a_7(x) &=& (1,0,0,0,0,0,0,0,-x_1,0,-x_2,-x_3)\ .
\end{eqnarray*}
The map\ \ \mb{b : (\Lambda^2Q)^{\oplus 2} \map (\mbb{C}^5)^*\otimes\Lambda^4Q = \Hom(\mbb{C}^5\otimes\mc{O}_{\mbb{P}_5},\Lambda^4Q)}\ \ in $(M)$ is defined by
\[
b(x)(\xi,\eta)(v) := -\eta\wedge(\Lambda^2\pi(x))(a_+(v)) + \xi\wedge(\Lambda^2\pi(x))(a_-(v))
\]
for\ \ \mb{v \in \mbb{C}^5}\ \ and\ \ \mb{\xi,\eta \in \Lambda^2Q(x)\ .}\ It is easily checked that the vectors
\begin{eqnarray*}
a_1(x) &=& (x_1x_2,0,x_1x_4,0,0,x_3x_4,0,0,0,0,0,0)\ ,\\
a_2(x) &=& (0,0,0,0,0,0,0,x_1x_3,0,x_2x_3,x_2x_4,0)
\end{eqnarray*}
are in $\ker(b(x))$, and that\ \ \mb{\underline{a} := (a_1,\ldots,a_7)}\ \ is a holomorphic frame field for $\ker(b)$ over\ \ \mb{U_0^* := \sett{x \in U_0}{x_i \not= 0\ ,\ 1 \leq i \leq 4}\ .}\ Since $(a_3,\ldots,a_7)$ is a basis of $\im(a)$, the projection\ \ \mb{\ker(b) \map \qmod{\ker(b)}{\im(a)} = E(-2)}\ \ maps $a_1,a_2$ to a holomorphic frame field\ \ \mb{\tilde{\underline{a}} := (\tilde{a}_1,\tilde{a}_2)}\ \ of $E(-2)$ over $U_0^*$. We write the matrix representation of $h_3(x)\vert_{\ker(b)}$ with respect to $\underline{a}$ as block matrix
\[
h_3(x)\vert_{\ker(b)} = \frac{1}{n}\left(
\begin{array}{cc}
C&\bar{B}^t\\B&A
\end{array}
\right)\ .
\]
where $A$ is the $5\times5$-matrix representing $h_3(x)\vert_{\im(a)}$. This can be calculated explicitely, using the matrix for $h_3(x)$; the result is
\[
C=
\left( 
\begin{array}{cc}
c_1 & 0 \\
0 & c_2 
\end{array}
\right)
\]
where 
\begin{eqnarray*}
c_1 &=& |x_1|^2|x_2|^2(1+|x_3|^2)+|x_1|^2|x_4|^2+|x_3|^2|x_4|^2(1+|x_2|^2)\ ,\\
c_2 &=& |x_1|^2|x_3|^2(1+|x_4|^2)+|x_2|^2|x_3|^2+|x_2|^2|x_4|^2(1+|x_1|^2)\ ,
\end{eqnarray*}
\[
B=
\left( 
\begin{array}{cc}
(|x_1|^2+|x_4|^2) \bar x_2 \bar x_3 & -(|x_2|^2+|x_3|^2) \bar x_1 \bar x_4 \\
(1+|x_2|^2)\bar x_3 \bar x_4 & -(|x_3|^2+|x_4|^2) \bar x_2 \\
(|x_1|^2+|x_3|^2)\bar x_4 & -(1+|x_4|^2) \bar x_1 \bar x_3 \\
(|x_2|^2+|x_4|^2)\bar x_1 & -(1+|x_1|^2) \bar x_2 \bar x_4 \\
(1+|x_3|^2) \bar x_1 \bar x_2 & -(|x_1|^2+|x_2|^2) \bar x_3  
\end{array}
\right)\ ,
\]
and
\[
A=
\left( 
\begin{array}{ccccc}
n+1 & \bar x_2 x_4 & \bar x_4 x_3 & \bar x_1 x_2 & \bar x_3 x_1 \\
\bar x_4 x_2 & n+|x_1|^2 & \bar x_3 & x_4 & \bar x_2 x_3 \\
\bar x_3 x_4 & x_3 & n+|x_2|^2 & \bar x_4 x_1 & \bar x_1 \\
\bar x_2 x_1 & \bar x_4 & \bar x_1 x_4 & n+|x_3|^2 & x_2 \\
\bar x_1 x_3 & \bar x_3 x_2 & x_1 & \bar x_2 & n+|x_4|^2
\end{array}
\right)\ .
\]
The matrix representation of the metric $h(x)$ in $E(x)$ with respect to $\tilde{\underline{a}}(x)$ is now given by
\[
h(x) = \frac{1}{n}(C - \bar{B}^t\cdot A^{-1}\cdot B)\ .
\]
We used MAPLE to explicitely calculate $h(x)$, but the resulting expression is to large to write down here.

We view\ \ \mb{x = (x_1,\ldots,x_4)}\ \ as holomorphic coordinates in $U_0$ via the standard chart\ \ \mb{(x_1,\ldots,x_4) \mapsto (1:x_1:\ldots:x_4)\ .}\ With respect to these coordinates, the K\"ahler form of the Fubini-Study metric $g$ is\ \ \mb{\omega_g = \frac{i}{2}\sum\limits_{\alpha,\beta = 1}^4 g_{\alpha\beta}dx_\alpha\wedge d\bar x_\beta\ ,}\ where\ \ \mb{g_{\alpha\beta} = \frac{\delta_{\alpha\beta}}{n} - \frac{\bar{x}_\alpha x_\beta}{n^2}}\ \ with\ \ \mb{n = 1 + \sum\limits_{i=1}^4\vert x_i\vert^2}\ \ as above.\\
Let $D$ be the Chern connection in $(E(-2),h)$, i.e. the unique $h$-unitary connection compatible with the holomorphic structure in $E(-2)$ (compare [K],[LT]), and\ \ \mb{F = D\circ D}\ \ its curvature. With respect to the holomorphic frame field $\tilde{\underline{a}}$ for $E(-2)$ over $U_0^*$, we write\ \ \mb{F = \left(F_{ij}\right)_{i,j = 1,2}}\ \ and\ \ \mb{F_{ij} = \sum\limits_{\alpha,\beta=1}^4 F_{ij\alpha\beta}dx_\alpha\wedge d\bar{x}_\beta\ ,}\ \mb{1 \leq i,j \leq 2\ .}\ Let be\ \ \mb{h = (H_{ij})_{i,j=1,2}}\ \ with respect to $\tilde{\underline{a}}$, and\ \ \mb{(H^{ij})_{i,j=1,2} := h^{-1}\ .}\ Then it holds 
\[
F_{ij\alpha\beta} = -\sum\limits_{k=1}^2 \frac{\partial^2 H_{ik}}{\partial x_\alpha\partial \bar{x}_\beta}H^{kj} +
\sum\limits_{k,l,m=1}^2 \frac{\partial H_{ik}}{\partial x_\alpha}H^{kl}\frac{\partial H_{lm}}{\partial \bar{x}_\beta}H^{mj}\ .
\]
The mean curvature $K$ of $h$ (with respect to $g$) is defined by the relation
\[
F\wedge\omega_g^3 = -\frac{i}{2}K\omega_g^4\ .
\]
With respect to $\tilde{\underline{a}}$ we write\ \ \mb{K = (K_{ij})_{i,j=1,2}\ ,}\ then it holds
\[
(*)\ \ \ K_{ij} = \sum\limits_{\alpha,\beta=1}^4 g^{\beta\alpha}F_{ij\alpha\beta} = -\sum\limits_{\alpha,\beta=1}^4 g^{\beta\alpha}\left(\sum\limits_{k=1}^2 \frac{\partial^2 H_{ik}}{\partial x_\alpha\partial \bar{x}_\beta}H^{kj} -
\sum\limits_{k,l,m=1}^2 \frac{\partial H_{ik}}{\partial x_\alpha}H^{kl}\frac{\partial H_{lm}}{\partial \bar{x}_\beta}H^{mj}\right)\ .
\]
where\ \ \mb{(g^{\alpha\beta})_{\alpha,\beta=1,\ldots,4} := \left((g_{\alpha\beta})_{\alpha,\beta=1,\ldots,4}\right)^{-1}\ ,}\ i.e.\ \ \mb{g^{\alpha\beta} = n(\delta_{\alpha\beta} + \bar{x}_\alpha x_\beta)\ .}

Since the $H_{ij}$ and $g^{\alpha\beta}$ are explicitely given, the calculation of $K(x^0)$ for a given point $x^0$ can now be done as follows (using MAPLE where necessary):\\
- determine $\frac{\partial^2 H_{ik}}{\partial x_\alpha\partial \bar{x}_\beta}(x)$, $\frac{\partial H_{ik}}{\partial x_\alpha}(x)$, $\frac{\partial H_{ik}}{\partial \bar{x}_\beta}(x)$, $1 \leq i,k \leq 2$, $1 \leq \alpha,\beta \leq 4$, for a general point $x$;\\
- substitute $x^0$ into $h$, and invert the scalar matrix $h(x^0)$ to get the $H^{ij}(x^0)$'s; \\
- substitute $x^0$ into $\frac{\partial^2 H_{ik}}{\partial x_\alpha\partial \bar{x}_\beta}$, $\frac{\partial H_{ik}}{\partial x_\alpha}$, $\frac{\partial H_{ik}}{\partial \bar{x}_\beta}$, $g^{\alpha\beta}$, $1 \leq i,k \leq 2$, $1 \leq \alpha,\beta \leq 4$; \\
- substitute the resulting scalars into the right hand side of equation $(*)$, and evaluate.

O.F.B. van Koert \\
M. L\"ubke

Mathematical Institute\\
Leiden University\\
PO Box 9512 \\
NL 2300 RA Leiden

okoert@math.leidenuniv.nl\\
lubke@math.leidenuniv.nl

\end{document}